\documentclass[11pt]{amsart}

\usepackage{amsmath, amssymb, amsthm, amsrefs}
\usepackage{latexsym}
\usepackage{indentfirst}
\usepackage{graphicx}
\usepackage{placeins}
\usepackage{booktabs}
\usepackage{algorithm}
\usepackage{algorithmic}
\usepackage{multirow}

\theoremstyle{plain}

\theoremstyle{definition}

\theoremstyle{remark}

\DeclareMathOperator{\diag}{diag}

\newcommand{\bbm}{\begin{bmatrix}}
\newcommand{\ebm}{\end{bmatrix}}

\newcommand{\R}{\mathrm{R}}

\newcommand{\p}{\partial}
\newcommand{\G}{\mathcal{G}}

\begin{document}

\title[Directional Preconditioner for High Frequency
  Scattering]{Directional Preconditioner for High Frequency Obstacle
  Scattering}

\author{Lexing Ying} 

\address{
  Department of Mathematics and Institute for Computational and Mathematical Engineering,
  Stanford University,
  Stanford, CA 94305
}

\email{lexing@math.stanford.edu}

\thanks{This work was partially supported by the National Science
  Foundation under award DMS-0846501 and the U.S. Department of
  Energy’s Advanced Scientific Computing Research program under award
  DE-FC02-13ER26134/DE-SC0009409. The author thanks Anil Damle for
  comments and suggestions.}



\begin{abstract}
  The boundary integral method is an efficient approach for solving
  time-harmonic obstacle scattering problems by a bounded
  scatterer. This paper presents the directional preconditioner for
  the iterative solution of linear systems of the boundary integral
  method. This new preconditioner builds a data-sparse approximation
  of the integral operator, transforms it into a sparse linear system,
  and computes an approximate inverse with efficient sparse and
  hierarchical linear algebra algorithms. This preconditioner is
  efficient and results in small and almost frequency-independent
  iteration counts when combined with standard iterative
  solvers. Numerical results are provided to demonstrate the
  effectiveness of the new preconditioner.
\end{abstract}

%\date{June 2014}

\maketitle

%-----------------------------------
\section{Introduction}

%continuous problem

This paper is concerned with solving the time-harmonic acoustic
obstacle scattering problems in two dimensions. Let $\Omega \subset
\R^2$ be a bounded scatterer with smooth boundary $\p\Omega$, $\omega$
be the frequency, and $e^{i\omega t} u_I(x)$ be the time-harmonic
incident wave. In the sound-soft scattering problem, the scattered
field $u(x)$ satisfies the Helmholtz equation with the Dirichlet
boundary condition
\begin{align*}
  \Delta u(x) + \omega^2 u(x) &= 0, && x\in \R^2 \setminus \Omega,\\
  u(x) &= -u_I(x), && x\in\p\Omega.
\end{align*}
In the sound-hard scattering problem, the scattered field $u(x)$
satisfies the Helmholtz equation with the Neumann boundary condition
\begin{align*}
  \Delta u(x) + \omega^2 u(x) &= 0, && x\in \R^2 \setminus \Omega,\\
  \frac{\partial u(x)}{\partial n(x)} &= -\frac{\partial u_I(x)}{\partial n(x)}, && x\in\p\Omega.
\end{align*}
In both cases, $u(x)$ satisfies the Sommerfeld radiation condition
\[
\lim_{r\rightarrow\infty} r^{1/2} \left(\frac{\p u}{\p r} - i \omega u\right) = 0.
\]

An effective way to solve these problems is the boundary integral
method and, more specifically, the combined field integral equation
(CFIE) \cite{colton-2013,nedelec-2001} formulation. This method relies
on the free space Green's function
\[
G(x,y) = \frac{i}{4} H^1_0(\omega |x-y|)
\]
of the Helmholtz equation. For the sound-soft scattering, we look for
a surface density $q(x), x\in \p\Omega$ such that for each
$x\in\p\Omega$,
\begin{equation}
  \frac{1}{2} q(x) + \int_{\p\Omega} \frac{\p G(x,y)}{\p n(y)} q(y) dy -
  i \eta \int_{\p\Omega} G(x,y) q(y) dy = -u_I(x).
  \label{eq:int1}
\end{equation}
Once $q(x)$ is computed, the scattered field $u(x)$ can be evaluated
through a boundary integral over $\p\Omega$. For the sound-hard
scattering, we look for $q(x), x\in \p\Omega$ such that for each
$x\in\p\Omega$,
\begin{equation}
  \frac{1}{2} q(x) -  \int_{\p\Omega} \frac{\p G(x,y)}{\p n(x)} q(y) dy +
  \frac{1}{i\eta} \oint_{\p\Omega} \frac{\p^2 G(x,y)}{\p n(x) \p n(y)} q(y) dy = 
  -\frac{\partial u_I(x)}{\partial n(x)}.
  \label{eq:int2}
\end{equation}
with $\oint$ being the principal value integral. In both cases, $\eta$
is typically chosen to be of order $O(\omega)$ (see \cite{kress-1985}
for example) and we refer to \cite{colton-2013,nedelec-2001} for
derivations and discussions of these integral equations.

By introducing the following operators
\begin{align*}
  (Sq)(x) &= \int_{\p\Omega} G(x,y) q(y) dy,\\
  (Dq)(x) &= \int_{\p\Omega} \frac{\p G(x,y)}{\p n(y)} q(y) dy,\\
  (D'q)(x)&= \int_{\p\Omega} \frac{\p G(x,y)}{\p n(x)} q(y) dy,\\
  (Nq)(x) &= \oint_{\p\Omega} \frac{\p^2 G(x,y)}{\p n(x) \p n(y)} q(y) dy,
\end{align*}
we can write \eqref{eq:int1} and \eqref{eq:int2} into operator forms:
\begin{align}
  &\left(\frac{1}{2} I + D - i\eta S\right) q = -u_I,\label{eq:dis1}\\
  &\left(\frac{1}{2} I + D' - \frac{1}{i\eta} N\right) q = -\frac{\partial u_I}{\partial n}.  \label{eq:dis2}
\end{align}

Standard approaches for discretizing these boundary integral equations
include the Nystr\"om method, the Galerkin method, and the collocation
method \cite{colton-2013,kress-2014}. To simplify the presentation, we
assume that the Nystr\"om method is used. For the other approaches,
the discussion remains similar as long as the basis functions employed
are local. A typical discretization of these integral equations
requires at least a couple of quadrature points per
wavelength. Assuming that both the diameter and the boundary length of
$\Omega$ are $\Theta(1)$, this implies that the boundary is
discretized with a set $P$ of $n = O(\omega)$ points. For the
resulting linear systems, we shall continue to use $S$, $D$, $D'$, and
$N$ to denote the discrete matrices associated with these
operators. Similarly, $q$, $u_I$, and $\p u_I/\p n$ are reused to
denote the discrete version of $q(x)$, $u_I(x)$, and $\p u_I(x)/\p
n(x)$ sampled at the quadrature points. Therefore, with this slight
abuse of notation, the discrete linear systems take the same form as
\eqref{eq:dis1} and \eqref{eq:dis2}.

%previous work
There has been a lot of work devoted to the fast solution of these
linear systems. Since the system is dense, the standard direct solvers
such as LU factorization take $O(n^3)$ steps, which is prohibitively
expensive. Recently, several linear-complexity approaches based on
recursive interpolative decomposition have been proposed by
\cite{corona-2013,martinsson-2005,ho-2013} for boundary integral
equations with non-oscillatory kernels. However, for high frequency
scattering where the kernel is oscillatory (i.e, $\omega=\Theta(n)$),
the complexity of these approaches is still cubic in $n$. The only
exception is for quasi-1D domains \cite{martinsson-2007} where the
boundary integral equation essentially reduces to the 1D case and the
complexity scales linearly in $\omega$.

For this reason, iterative methods such as GMRES and TFQMR
\cite{freund-1993,saad-1986,saad-2003} are the main approaches for
solving these problems. In these cases, though the CIFEs have much
better conditioning properties compared to other integral
formulations, the number of iterations can grow quickly with $\omega$.
Therefore, for high frequency scattering problems, there is a clear
need for improving the conditioning properties of these operators.

Over the past twenty years, there has been a significant amount of
research devoted to this task.  A couple of algorithms suggest
improving the conditioning property via modifying the standard CFIE
formulation. For example, one line of work is to replace the $i\eta$
term in \eqref{eq:dis1} with better approximations of the
Dirichlet-to-Neumann (DtN) operator and the $1/(i\eta)$ term in
\eqref{eq:dis2} with better approximations of the Neumann-to-Dirichlet
(NtD) operator
\cite{alouges-2007,antoine-2005,antoine-2007,bruno-2012}. Typically,
these new approximations are derived from leading order terms of the
pseudo-differential symbols of the DtN and NtD operators.

A second approach is to precondition the integral equation. Most work
here considers the electric field integral equation for
electromagnetic scattering and follows the famous Calderon
relationship \cite{steinbach-1998,christiansen-2002,antoine-2004}. The
resulting integral equations are of Fredholm second kind with good
conditioning properties. However, the number of matrix vector
multiplications per iteration is doubled.
%However, there is little work for preconditioning the CIFEs.

There has also been a lot of work on sparsifying the integral
operators using special basis functions, such as local cosine bases
\cite{averbuch,bradie} and optimized wavelet packets
\cite{deng1,deng2,golik,huybrechs}. The resulting sparse
representations typically have $O(n^{4/3})$ non-zero entries. Recently
in \cite{demanet-2010}, an approach using the wave atom transform
\cite{demanet-2007} results a sparse representation with $O(n\log n)$
non-zero entries. In \cite{canning}, Canning claimed to obtain a
sparse approximation with $O(n)$ non-zero entries via locally
mollified exponential functions. However, when good accuracy is
required, most of these methods access all entries of the integral
operator, thus requiring an $O(n^2)$ precomputation cost to assemble
the whole matrix.

In this paper, we propose a new method for preconditioning the CFIEs
by incorporating the ideas from sparse representation. This approach
builds a data-sparse representation of the boundary integral operator,
transforms it into a sparse linear system, and computes an approximate
inverse with efficient sparse and hierarchical linear algebra
algorithms. This preconditioner is highly efficient to construct and
to apply. It results in small and almost frequency-independent
iteration counts when combined with standard iterative solvers. The
rest of the paper is organized as follows. Section 2 describes the
algorithm and Section 3 presents the results. Future work and open
questions are discussed in Section 4.

%-----------------------------------
\section{Algorithm}

For frequency $\omega$, the wavelength $\lambda$ is $2\pi/\omega$. We
assume that the scatterer boundary $\p\Omega$ is $C^2$ and both the
diameter and the boundary length of $\Omega$ is $\Theta(1)$.  To
simplify the discussion, we suppose that the length $L$ of $\p\Omega$
is equal to $4^q \lambda$ where $q$ is a positive integer. The actual
number $4^q$ is not essential but it makes the presentation
simpler. Combining this with $L=\Theta(1)$ implies that $\omega =
O(4^q)$.

Suppose that $\rho:\partial\Omega \rightarrow [0,L]$ is the arclength
parametrization of the boundary and that the boundary is sampled with
$n=4^q p$ discretization points for some $p=\Theta(1)$, i.e.,
$p$ points per wavelength $\lambda$.

%---------
\subsection{Data-sparse approximation}

We start by decomposing the boundary into sufficiently planar
segments. Initially, the boundary is partitioned into $2^q$ segments,
each of length $2^{q}\lambda$ and with $2^{q}p$ points. Each such
segment is further partitioned hierarchically until one of two
situations happen:
\begin{itemize}
\item First, it is stopped if the length of the segment is bounded by
  $2^{q}\lambda/\sqrt{c}$ where $c$ is the maximum absolute value of
  the curvature in the current segment. Such a segment is called {\em
    almost-planar}.
\item Second, it is stopped when the length of the segment is bounded
  by $m_\ell \lambda$. Typically $m_\ell=2$ or $4$. Such a segment is
  called a non-planar leaf.
\end{itemize}
We denote the final set of segments by $\G = \{P_1,\ldots,P_m\}$,
where the segments $P_i$ are ordered according to their positions on
the boundary. Notice that since the boundary is assumed to be $C^2$,
all segments in $\G$ are almost-planar for sufficiently large
$\omega$. Therefore, in the following discussion, it is safe to regard
all segments $P_i$ as almost-planar.

The discussion here shall treat the sound-soft case \eqref{eq:dis1}
and the sound-hard case \eqref{eq:dis2} in the same way, since the
kernels of these two have the same oscillatory pattern. Therefore, it
is convenient to use the general form
\begin{equation}
M q = f
\label{eq:Mqf}
\end{equation}
for both of them in the discussion. Based on how the segments are
generated, each $P_j$ is of length $2^{\ell_j} \lambda$ for some
integer $\ell_j$ and contains $2^{\ell_j} p$ equally-spaced
discretization points. After ordering the unknowns according to the
ordering of $P_j$, the matrix $M$ can be written as the following
block form
\[
M = \begin{bmatrix}
  M_{11} & \ldots & M_{1m} \\
  \vdots & \ddots & \vdots \\
  M_{m1} & \ldots & M_{mm}
\end{bmatrix},
\]
where $M_{ij}$ is of size $2^{\ell_i}p \times 2^{\ell_j}p$.  The next
step is find a data-sparse approximation for the blocks $M_{ij}$.

\subsubsection{Diagonal blocks}
Let us first consider a diagonal block $M_{jj}$, which represents the
interaction between $P_j$ and itself. Since $P_j$ is almost planar, we
can treat it approximately as flat. Therefore, we have 
\[
M_{jj} \approx B_j,
\]
where $B_j$ is obtained by restricting the integral operator to a
straight segment of length $2^{\ell_j}\lambda$ with $2^{\ell_j}p$
equally spaced quadrature points. Noticing that $B_j$ only depends on
$\ell_j$ and that there are only a few choices for $\ell_j$, we can
clearly precompute these matrices. Going through all $P_i$ gives the
following approximation $B$ to the block-diagonal part of $M$,
\[
B = \begin{bmatrix}
  B_1 & & \\
  & \ddots & \\
  & & B_m\\
\end{bmatrix},
\]

\subsubsection{Off-diagonal blocks}
\begin{figure}[h!]
  \begin{center}
    \includegraphics[height=1in]{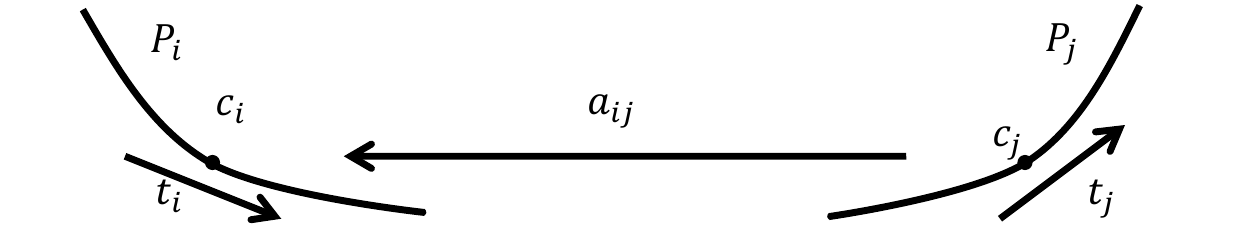}    %\hspace{0.25in}
  \end{center}
  \caption{The geometric setup of two segments $P_i$ and $P_j$.}
  \label{fig:lr}
\end{figure}

Next we consider the off-diagonal blocks, i.e, $M_{ij}$ with
$i\not=j$. We define
\begin{itemize}
\item $c_i$ and $c_j$ to be the centers of segments $P_i$ and $P_j$,
\item $t_i$ and $t_j$ to be the tangent directions of $\p\Omega$ at
  the centers of $P_i$ and $P_j$, and
\item $a_{ij}$ to be the unit direction from $c_j$ to $c_i$, i.e.,
  $(c_i-c_j)/|c_i-c_j|$ (see Figure \ref{fig:lr} for an illustration).
\end{itemize}

The main difficulty of working with high frequency Helmholtz kernel
$G(x,y)$ is its highly oscillatory behavior. Since the oscillation of
$G(x,y)$ and its derivatives come from the term $\exp(i\omega|x-y|)$,
it is instructive to focus on $\exp(i\omega|x-y|)$ for a moment. In
the following discussion, the sign $\sim$ is used to denote an
approximation up to a non-oscillatory multiplicative term.

A first observation is that 
\begin{equation}
  \exp(i\omega|x-y|) \sim \exp(i\omega a_{ij}\cdot (x-y)).
  \label{eq:stt}
\end{equation}
To see this, assume without loss of generality that $a_{ij}=(1,0)^t$
by rotating the coordinates accordingly and placing the origin on
$\p\Omega$ between $P_i$ and $P_j$. In the new coordinate system, we
have
\begin{align*}
  \exp(i\omega|x-y|) & = \exp(i\omega(x_1-y_1)) \exp(i\omega(|x-y|-(x_1-y_1)))\\
  & = 
  \exp(i\omega(x_1-y_1))
  \exp\left(i\omega(x_1-y_1) \left(\sqrt{1+\frac{|x_2-y_2|^2}{|x_1-y_1|^2}}-1\right)\right),
\end{align*}
where $x=(x_1,x_2)^t$ and $y=(y_1,y_2)^t$. Up to a constant factor,
the phase of the last term can be estimated with $\omega
\frac{|x_2-y_2|^2}{|x_1-y_1|}$.
\begin{itemize}
\item When the segments $P_i$ and $P_j$ are $\Theta(1)$ distance from each
  other, we estimate this by
  \[
  \omega \frac{|x_2-y_2|^2}{|x_1-y_1|} \lesssim \omega (2^q\lambda)^2 = \Theta(1).
  \]
\item When $P_i$ and $P_j$ are close to each other, we use quadratic
  approximation in the rotated frame $|x_2| \lesssim |x_1|^2$ and
  $|y_2| \lesssim |y_1|^2$ and the fact that $x_1$ and $y_1$ has
  different signs to conclude
  \[
  \omega \frac{|x_2-y_2|^2}{|x_1-y_1|} \lesssim \omega |x_1-y_1|^3 \lesssim \omega
  (2^q\lambda)^3 = o(1).
  \]
\end{itemize}
In both cases, $\omega \frac{|x_2-y_2|^2}{|x_1-y_1|}$ is bounded and
hence \eqref{eq:stt} is valid up to a non-oscillatory multiplicative
term. Next, we rewrite
\begin{align}
  & \exp(i\omega|x-y|) \sim \exp(i\omega a_{ij}\cdot(x-y)) \nonumber \\
  &=\exp(i\omega a_{ij}\cdot((x-c_i)+(c_i-c_j)+(c_j-y)) ) \nonumber \\
  &=\exp(i\omega a_{ij}\cdot(x-c_i)) \cdot \exp(i\omega a_{ij}\cdot (c_i-c_j)) 
  \cdot \exp(-i\omega a_{ij}\cdot(y-c_j)). \label{eq:expsep}
\end{align}

To approximate the first term in \eqref{eq:expsep}, we perform a Taylor
expansion for $\rho^{-1}(t)$ near $t=\rho(c_i)$ and evaluate it at
$\rho(x)$:
\begin{align*}
  & \left| \rho^{-1}(\rho(x)) - \left( \rho^{-1}(\rho(c_i)) + t_i (\rho(x)-\rho(c_i)) \right) \right|\\
  & \lesssim \frac{1}{2} |\rho(x)-\rho(c_i)|^2 c \le \frac{1}{2} (2^\ell/\sqrt{c}\lambda)^2 c = O(\lambda),
\end{align*}
where $c$ is the maximum absolute value of the curvature in $T$. The
inequality here uses the fact that $P_i$ is almost-planar.  This is
equivalent to
\[
(x-c_i) = (\rho(x)-\rho(c_i))\cdot t_i + O(\lambda).
\]
Multiplying it with $i \omega a_{ij}$ and taking exponential gives the approximation
\begin{equation}
  \exp( i\omega a_{ij} \cdot (x-c_i)) \sim \exp(i (\omega a_{ij}\cdot t_i)(\rho(x)-\rho(c_i)))
  \label{eq:iflat}
\end{equation}
Since $S$ is almost-planar, the same argument works for $(y-c_j)$ and
gives
\[
(y-c_j) = (\rho(y)-\rho(c_j))\cdot t_j + O(\lambda)
\]
and
\begin{equation}
  \exp(-i\omega a_{ij}\cdot (y-c_j)) \sim \exp(i (-\omega a_{ij}\cdot t_j)(\rho(y)-\rho(c_j))).
  \label{eq:jflat}
\end{equation}

We can now further approximate the phase function of the complex
exponentials in \eqref{eq:iflat} and \eqref{eq:jflat} as
follows. Noticing that $\omega a_{ij}\cdot t_i \in [-\omega,\omega]$,
we partition the interval $[-\omega,\omega]$ into $2^{\ell_i+1}$
equally spaced subintervals with a set $K_i$ of $2^{\ell_i+1}+1$
gridpoints. We define $[k]_i$ to the value of rounding $k$ to the
nearest gridpoint in $K_i$. Then
\[
(\omega a_{ij}\cdot t_i - [\omega a_{ij}\cdot t_i]_i)(\rho(x)-\rho(c_i)) \le
\frac{2\omega}{2^{\ell_i+1}} \cdot \frac{1}{2}\cdot \frac{2^{\ell_i}\lambda}{2} = \frac{2\pi}{4} = O(1).
\]
Thus, replacing the phase $\omega a_{ij}\cdot t_i$ with $[\omega
  a_{ij}\cdot t_i]_i$ in \eqref{eq:iflat} introduces an extra
non-oscillatory term
\begin{equation}
  \exp( i\omega a_{ij} \cdot (x-c_i)) \sim \exp(i [\omega a_{ij}\cdot t_i]_i (\rho(x)-\rho(c_i)))
  \label{eq:iflatnew}
\end{equation}

Similarly, we partition the interval $[-\omega,\omega]$ into
$2^{\ell_j}$ equal pieces with a set $K_j$ of $2^{\ell_j+1}+1$
gridpoints. By defining $[k]_j$ to the value of rounding $k$ to the
nearest gridpoint in $K_j$, we again have for $y\in S$
\[
(-\omega a_{ij}\cdot t_j - [-\omega a_{ij}\cdot t_j]_j)(\rho(y)-\rho(c_j)) \le
\frac{2\omega}{2^{\ell_j+1}} \cdot \frac{1}{2}\cdot \frac{2^{\ell_j}\lambda}{2} = \frac{2\pi}{4} = O(1),
\]
This change of the phase function also introduces an extra
non-oscillatory term
\begin{equation}
  \exp(-i\omega a_{ij}\cdot (y-c_j)) \sim \exp(i [-\omega a_{ij}\cdot t_j]_j (\rho(y)-\rho(c_j))).
  \label{eq:jflatnew}
\end{equation}

By introducing
\begin{align}
  & k^i_{ij} = [\omega a_{ij} t_i]_i, \nonumber \\
  & k^j_{ij} = [-\omega a_{ij} t_j]_j,\nonumber\\
  & U_i(x,k) = \exp(i k (\rho(x)-\rho(c_i))),\nonumber\\
  & U_j(y,k) = \exp(i k (\rho(y)-\rho(c_j))),\label{eq:Ui}
\end{align}
and putting \eqref{eq:iflatnew} and \eqref{eq:jflatnew} in
\eqref{eq:expsep}, we have the following approximation
\[
\exp(i\omega|x-y|) \sim U_i(x,k^i_{ij}) \cdot 
\exp(i\omega (c_i-c_j) a_{ij}) \cdot U_j(y,k^j_{ij}).
\]
for all $x\in P_i$ and $y\in P_j$.

Since the kernel $G(x,y)$ and its derivatives have the same
oscillation pattern as $\exp(i\omega|x-y|)$, $U_i(x,k^i_{ij})$ and
$U_j(y,k^j_{ij})$ also capture the oscillations of $G(x,y)$ for $x\in
P_i$ and $y\in P_j$.  Therefore, repeating the same argument gives the
following representation of the block $M_{ij}$:
\begin{equation}
M_{ij}(x,y) = U_i(x,k^i_{ij}) \cdot \tilde{M}_{ij}(x,y) \cdot
U_j(y,k^j_{ij})
\label{eq:Mtilde}
\end{equation}
for $x\in P_i$ and $y\in P_j$, where the non-oscillatory term
$\tilde{M}_{ij}(x,y)$ is defined through this representation. Since
$\tilde{M}_{ij}(x,y)$ is non-oscillatory, we can approximate it with
Chebyshev interpolation. For this, we define
\begin{itemize}
\item $R_i$ and $R_j$ to be the Chebyshev grids of a constant size
  $m_c$ in $P_i$ and $P_j$, respectively, and
\item $I_i$ and $I_j$ to be the corresponding interpolation operators,
  with entries given by $I_i(x,b)$ for $b\in R_i$ and $I_j(y,b)$ for
  $b\in R_j$.
\end{itemize}
This results the following approximation
\[
\tilde{M}_{ij}(x,y) \approx I_i \cdot \tilde{M}_{ij}(R_i,R_j) \cdot I_j^t.
\]
Putting this together with \eqref{eq:Mtilde} and using matrix form
gives
\[
M_{ij} \approx \diag(U_i(:,k^i_{ij})) \cdot I_i \cdot \tilde{M}_{ij}(R_i,R_j) \cdot I_j^t \cdot 
\diag(U_j(:,k^j_{ij})).
\]

For the data-sparse representation and the preconditioner, we need an
aggressive rank-1 approximation for $M_{ij}$ of form
\[
M_{ij} \approx U_i(:,k^i_{ij})) e_{ij} U_j(\cdot,k^j_{ij})^t
= \diag(U_i(:,k^i_{ij})) \cdot w_i \cdot e_{ij} \cdot w_j^t \cdot \diag(U_j(\cdot,k^j_{ij})),
\]
where $e_{ij}$ is a constant to be determined and $w_i$ and $w_j$ are
the all-one vectors of length $2^{\ell_i}p$, and $2^{\ell_j}p$,
respectively. To determine $e_{ij}$, we can solve for it from a least
square problem
\begin{equation}
  e_{ij} = \text{argmin}_{e} 
  \|I_i \cdot \tilde{M}_{ij}(R_i,R_j) \cdot I_j^t - w_i \cdot e \cdot w_j^t \|^2.
  \label{eq:entry}
\end{equation}
The solution is 
\[
e_{ij} = (w_i^\dagger \cdot I_i) \cdot \tilde{M}_{ij}(R_i,R_j) \cdot (I_j^t \cdot (w_j^t)^\dagger).
\]
Notice that $(w_i^\dagger \cdot I_i)$ and $(I_j^t \cdot
(w_j^t)^\dagger)$ only depend on $\ell_i$ and $\ell_j$
respectively. Therefore, they can be precomputed and the remaining
cost of computing $e_{ij}$ is equal to $O(m_c^2)$.

Going through all pairs $(i,j)$ with $i\not=j$ yields the following
approximation for the off-diagonal part of $M$:
\[
U E U^t.
\]
Here
\[
U = \begin{bmatrix}
  U_1 & & \\
  & \ddots & \\
  & & U_m
\end{bmatrix},
\]
where $U_i$ is a matrix of size $2^{\ell_i}p \times (2^{\ell_i+1}+1)$
given by \eqref{eq:Ui}. The $E$ matrix also has a $m\times m$ block
form
\[
E = \begin{bmatrix}
E_{11} & \ldots & E_{1m}\\
\vdots & \ddots & \vdots\\
E_{m1} & \ldots & E_{mm}\\
\end{bmatrix}
\]
where $E_{ij}$ is a matrix of size $(2^{\ell_i+1}+1) \times
(2^{\ell_j+1}+1)$ with rows and columns indexed by $K_i$ and $K_j$.
$E_{ij}$ is a matrix with value $e_{ij}$ at entry
$(k^i_{ij},k^j_{ij})$ and zero everywhere else. Here we emphasize that
\begin{itemize}
\item $U_i$ is a partial Fourier matrix, and 
\item $E$ is extremely sparse.
\end{itemize}
These observations turn out to be essential in the construction of the
preconditioner.

Summarizing the discussion for both the diagonal and off-diagonal
blocks, we hold the data-sparse approximation
\begin{equation}
M \approx B + U E U^t.
\label{eq:Mapp}
\end{equation}

%---------
\subsection{Directional preconditioner}

To precondition \eqref{eq:Mqf}, we use the approximation
\eqref{eq:Mapp} and consider the solution $q$ of
\[
(B + UEU^t) q = f.
\]
First, introducing new vectors $r=-U^tq$ and $p=-Er$ gives an
equivalent augmented system
\begin{equation}
\begin{bmatrix}
  B   & U & 0\\
  U^t & 0 & I\\
  0 & I & E
\end{bmatrix}
\begin{bmatrix}
  q \\p \\r
\end{bmatrix}
=
\begin{bmatrix}
  f \\0 \\0
\end{bmatrix}.
\label{eq:aug}
\end{equation}
Factorizing the matrix in \eqref{eq:aug} gives
\begin{equation}
\begin{bmatrix}
  I   &  & \\
  U^tB^{-1} & I & \\
  & & I
\end{bmatrix}
\begin{bmatrix}
  I   &  & \\
  & I & \\
  & -T & I
\end{bmatrix}
\begin{bmatrix}
  B &    &\\
    & -S & \\
    &    & W
\end{bmatrix}
\begin{bmatrix}
  I   &  & \\
  & I & -T\\
  & & I
\end{bmatrix}
\begin{bmatrix}
  I & B^{-1}U & \\
  & I & \\
  & & I
\end{bmatrix},
\label{eq:fac}
\end{equation}
with
\[
S = U^t B^{-1} U,\quad
T = S^{-1},\quad
W = E+T.
\]
Since both $U$ and $B$ are block-diagonal, $S$ and $T$ are also
block-diagonal
\[
S = \begin{bmatrix}
  S_1 & & \\
  & \ddots & \\
  & & S_m
\end{bmatrix},\quad
T = \begin{bmatrix}
  T_1 & & \\
  & \ddots & \\
  & & T_m
\end{bmatrix},
\]
with $S_i = U_i^t B_i^{-1} U_i$ and $T_i = S_i^{-1}$. Inverting the
factorization \eqref{eq:fac} gives
\begin{align}
  \begin{bmatrix}
    q \\ p \\r
  \end{bmatrix}
  = &
  \begin{bmatrix}
  I   & -B^{-1}U & \\
  & I & \\
  & & I
  \end{bmatrix}
  \begin{bmatrix}
    I   &  & \\
    & I & T\\
    & & I
  \end{bmatrix}
  \begin{bmatrix}
    B^{-1}   & &\\
    & -T & \\
    &  & W^{-1}
  \end{bmatrix} \nonumber\\
  &
  \begin{bmatrix}
    I   &  & \\
    & I & \\
  & T & I
  \end{bmatrix}
  \begin{bmatrix}
    I & & \\
    -U^tB^{-1} & I & \\
    & & I
  \end{bmatrix}
  \begin{bmatrix}
    f \\ 0 \\0
  \end{bmatrix}.
  \label{eq:esol}
\end{align}

\newcommand{\app}[1]{\lfloor#1\rfloor}

%Since $B$, $B^{-1}$, $S$, and $T$ are all block diagonal matrices,
%for them we focus on the diagonal blocks.

\begin{figure}[h!]
  \begin{center}
    \includegraphics[height=1.5in]{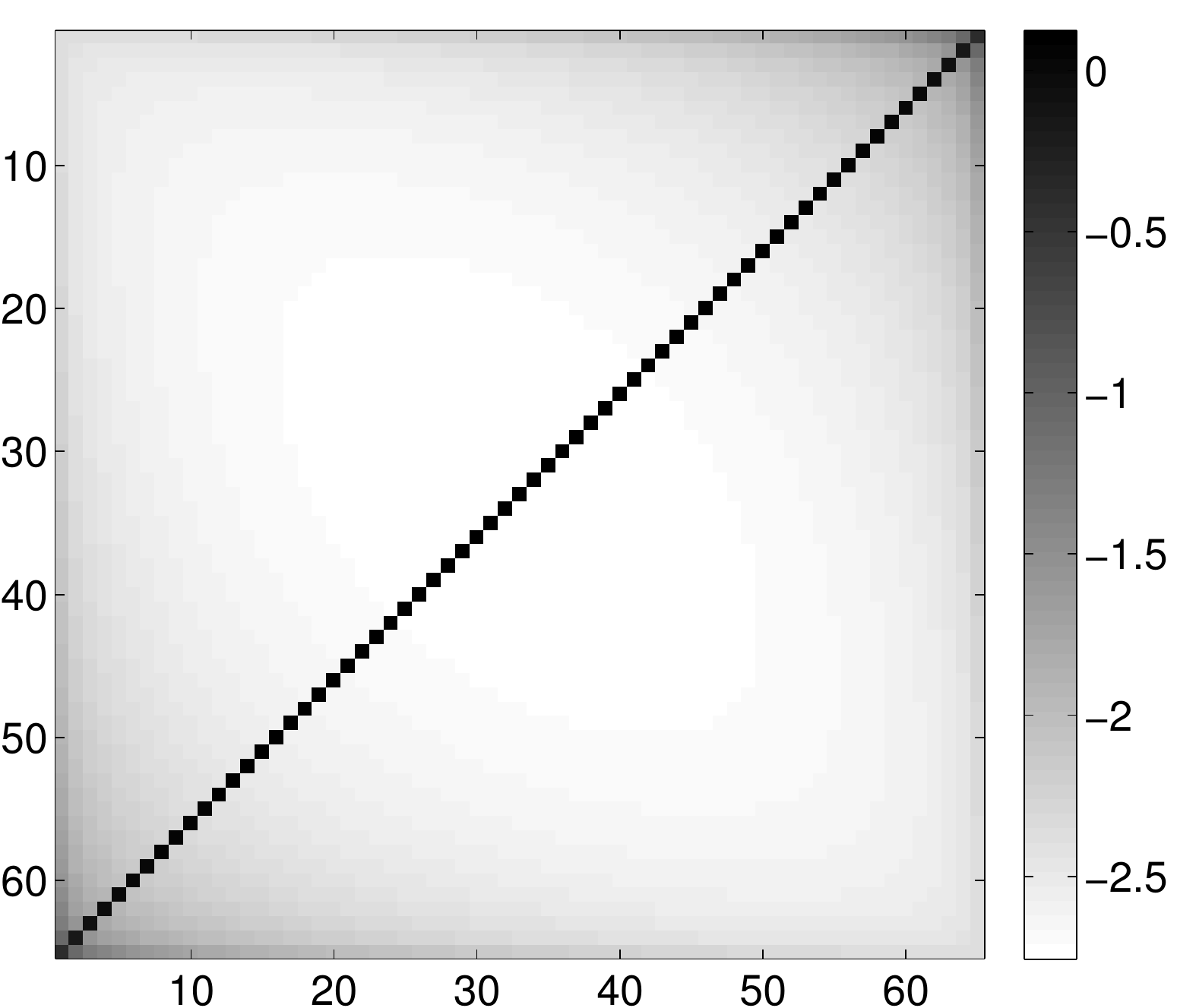}    \hspace{0.25in}    \includegraphics[height=1.5in]{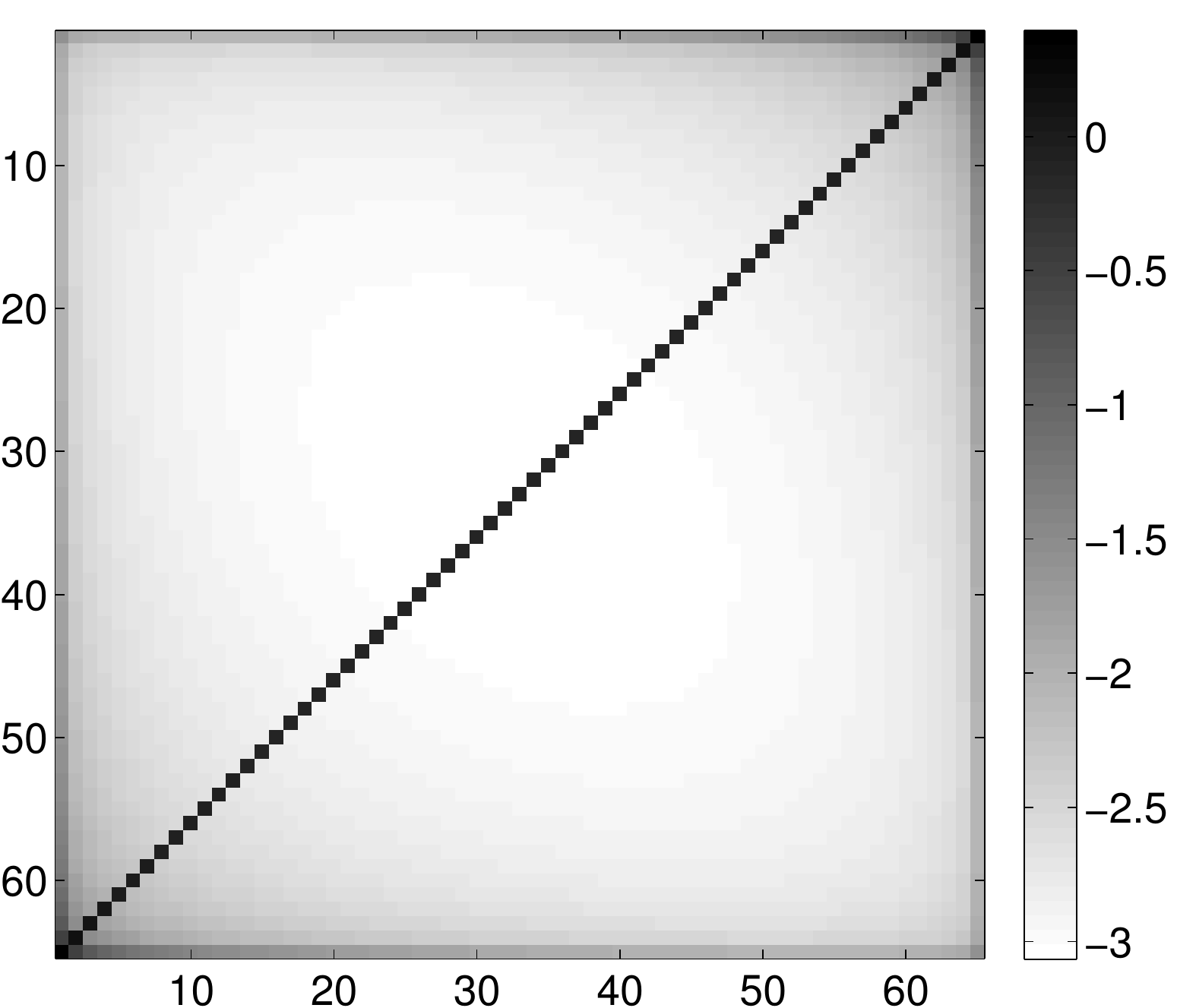}    \\
    \includegraphics[height=1.5in]{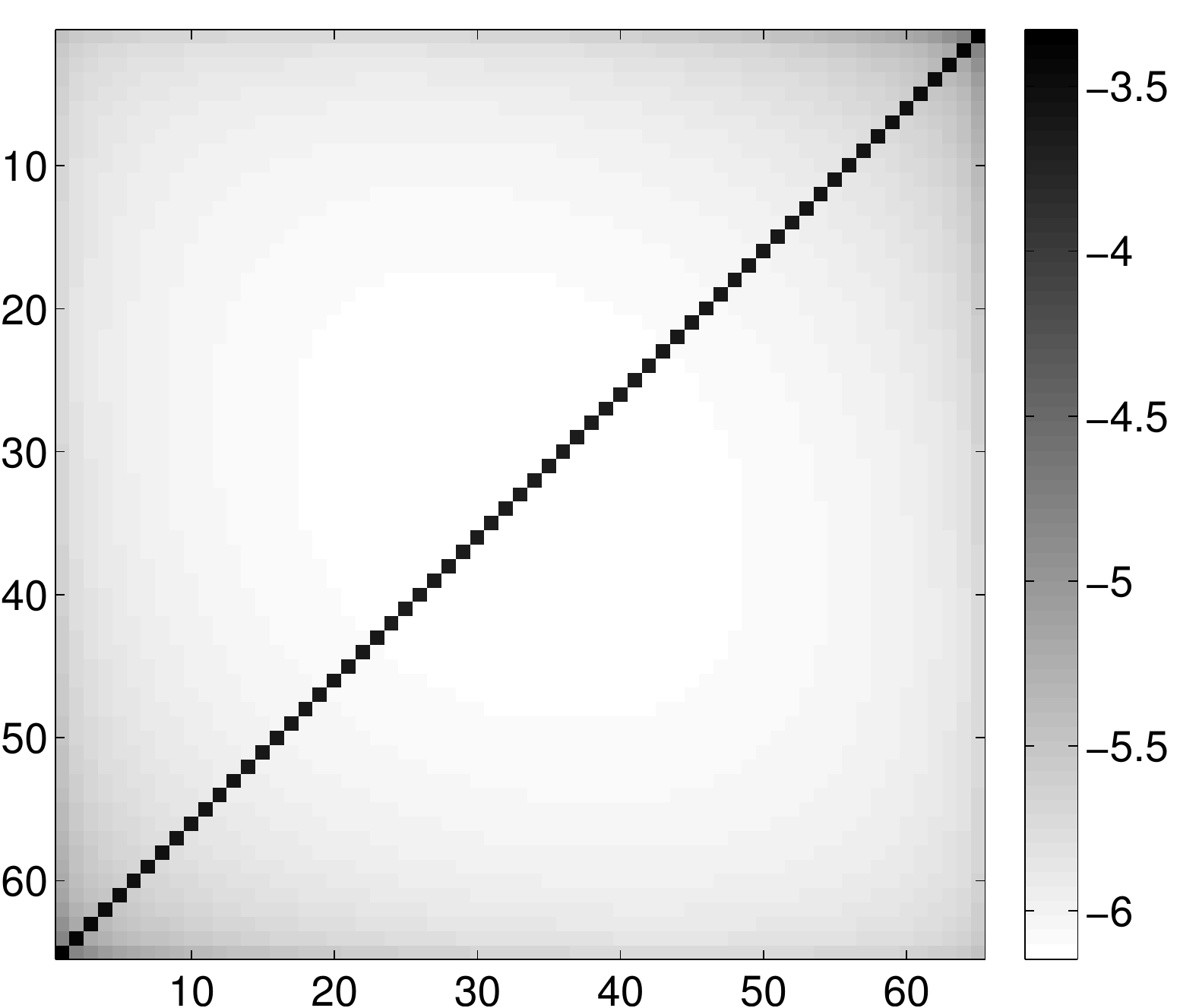}    \hspace{0.25in}    \includegraphics[height=1.5in]{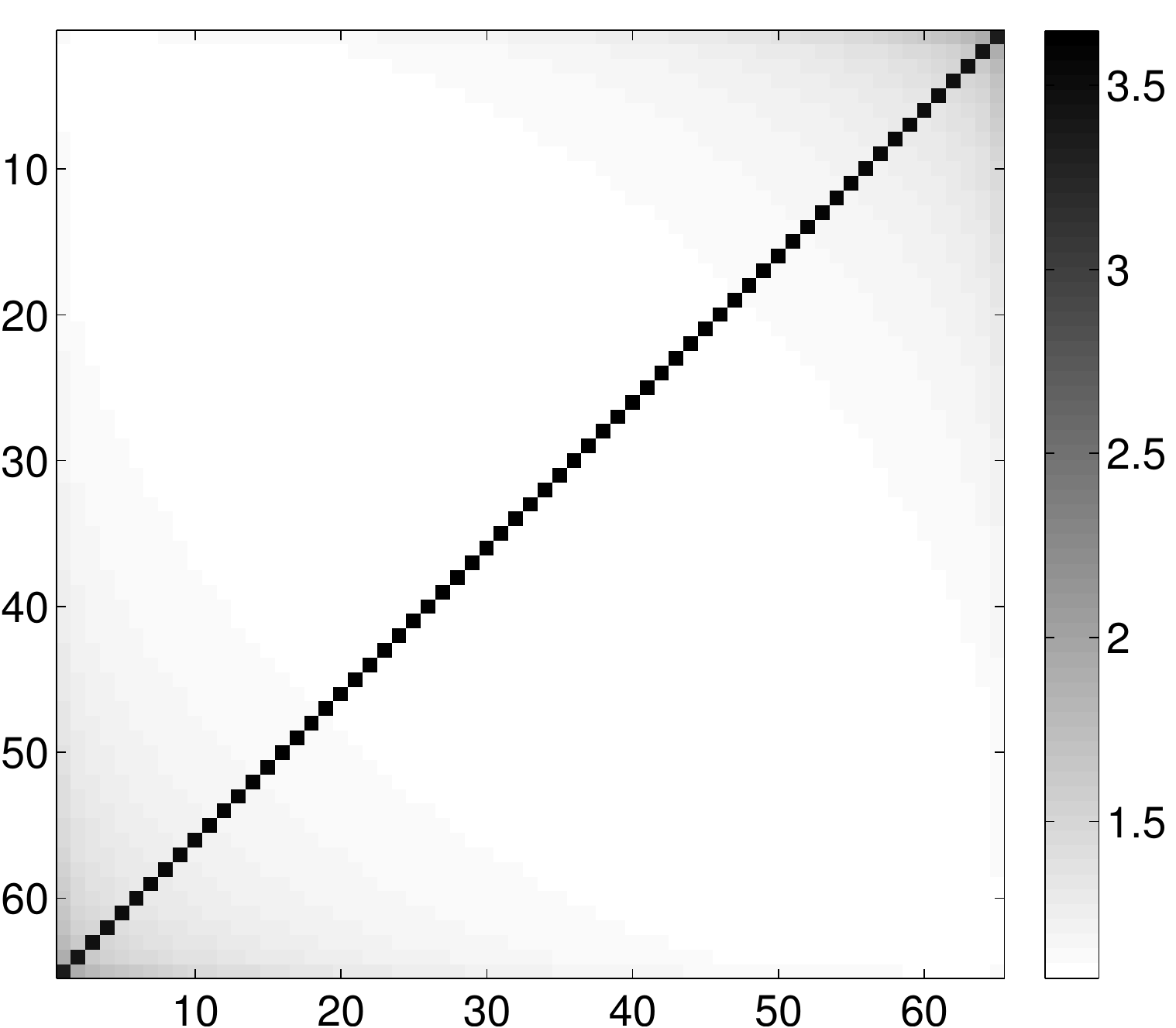}    \\
  \end{center}
  \caption{The absolute value of the entries of $S_i$ (left) and $T_i$
    (right) in logarithmic scale (base $10$) for the sound-soft case
    (top row) and the sound-hard case (bottom row). For this specific
    segment, $\ell_i = 5$ and $p=8$. Hence, there are 256
    equally-spaced quadrature points in $P_i$ and the cardinality of
    $K_i$ is 65. Both $S_i$ and $T_i$ are numerical sparse and
    anti-diagonally dominant.}
  \label{fig:ST}
\end{figure}

Applying \eqref{eq:esol} exactly can be quite costly. In order to
construct an efficient preconditioner, it is essential to approximate
\eqref{eq:esol} aggressively while without sacrificing too much
accuracy. For a matrix $A$, we shall use the notation
$\app{A}$ to stand for its approximation, but the actual approximation
scheme can be different for different matrices.
\begin{itemize}
\item First, the $2^{\ell_i}p\times 2^{\ell_i}p$ matrix $B_i$ is the
  restriction of the integral operator to a straight segment of length
  $2^{\ell_i}\lambda$. Since the geometry is fixed and there are only
  a few choices for $\ell_i$, all $B_i$ and $B_i^{-1}$ can be
  precomputed. Since this is also a one-dimensional problem (i.e.,
  restriction to a straight segment), we can use the hierarchical
  matrix algebra \cite{borm-2006} or the hierarchical semi-separable
  (HSS) matrices \cite{xia-2010} to compress and apply $B_i^{-1}$
  efficiently. We denote the approximation of $B_i^{-1}$ with
  $\app{B_i^{-1}}$ and accordingly $\app{B^{-1}}$ for $B^{-1}$.
\item Second, since each $U_i$ is a partial Fourier matrix, applying
  $U$ and $U^t$ reduces to a number of FFTs, which is highly
  efficient.
\item Third, $S_i$ and $T_i$ can be precomputed as they only depend on
  $U_i$ and $B_i$, both of which have already been precomputed. An
  important observation is that $T_i$ is numerically sparse (see
  Figure \ref{fig:ST}). Therefore, for the sake of efficiency, we
  approximate $T_i$ with $\app{T_i}$, which is obtained by
  thresholding the entries in absolute value. In the numerical
  results, the number of non-zero entries in $\app{T_i}$ is kept
  proportional to the dimension of $T_i$.  This approximation of $T$
  is denoted by $\app{T}$.
\item 
  The final task is to build an approximate inverse of $W=E +
  T$. Here, the essential observation is that
  \begin{itemize}
  \item $T$ concentrates on its anti-diagonal (see Figure
    \ref{fig:ST}).
  \end{itemize}
  Based on this, we define $\app{W}$ to be the sum of $E$ and the
  anti-diagonal of $T$ (i.e., thresholding the rest entries of $T$ to
  zero). $\app{W}$ is extremely sparse as the number of non-zeros is
  about $3/2$ times the dimension of the matrix in most cases. Hence,
  we perform a sparse LU decomposition and set
  \[
  \app{W} = L_{\app{W}} R_{\app{W}},
  \]
  where $L_{\app{W}}$ and $R_{\app{W}}$ are sparse upper and lower
  triangular matrices up to possible permutations.
\end{itemize}

Once these approximations have been formed, we approximate \eqref{eq:esol}
with
\begin{align}
\begin{bmatrix}
q \\ p \\r
\end{bmatrix}
\Leftarrow &
\begin{bmatrix}
  I & -\app{B^{-1}}U & \\
  & I & \\
  & & I
\end{bmatrix}
\begin{bmatrix}
  I   &  & \\
  & I & \app{T}\\
  & & I
\end{bmatrix}
\begin{bmatrix}
  \app{B^{-1}}   & &\\
  & -\app{T} & \\
  &  & R_{\app{W}}^{-1}L_{\app{W}}^{-1}
\end{bmatrix}\nonumber \\
&\begin{bmatrix}
  I & & \\
  & I & \\
  & \app{T} & I
\end{bmatrix}
\begin{bmatrix}
  I & & \\
  -U^t\app{B^{-1}} & I & \\
  & & I
\end{bmatrix}
\begin{bmatrix}
f \\ 0 \\0
\end{bmatrix}.
\label{eq:eapp}
\end{align}
%where $\app{R}^{-1}$ and $\app{L}^{-1}$ are applied through sparse
%backward and forward substitutions.
We emphasize again that the following approximations are used for
computing \eqref{eq:eapp}:
\begin{itemize}
\item replacing $B^{-1}$ with $\app{B^{-1}}$ via hierarchical matrix
  or HSS approximation for each $B_i^{-1}$,
\item applying $U$ rapidly via fast Fourier transform for each $U_i$,
\item replacing $T$ with $\app{T}$ via sparse approximation for each
  $T_i$, and
\item replacing $W^{-1}$ with $R_{\app{W}}^{-1}L_{\app{W}}^{-1}$ via
  sparse backward and forward substitutions for $L_{\app{W}}$ and
  $R_{\app{W}}$.
\end{itemize}

Based on \eqref{eq:eapp}, our preconditioner is defined as
follows. For a given $f$, it
\begin{itemize}
\item forms vector $(f^t,0,0)^t$,
\item carries out the computation of \eqref{eq:eapp}, and
\item extracts the first component $q$ of the resulting vector.
\end{itemize}
Since the key step of constructing a data-sparse representation of the
operator relies on the directional nature of the kernel $G(x,y)$, we
name it {\em directional preconditioner}.

%---------
\subsection{Complexity analysis}

We first consider the setup cost the preconditioner, i.e., the
approximate factorization in \eqref{eq:eapp}. Since the matrices
$B_i^{-1}$ and $T_i$ for a segment $P_i$ only depend on the integer
length parameter $\ell_i$ of $P_i$, the possible choices for these
matrices are fixed and independent of the scatterer. All these
possible choices can be precomputed once and for all and stored for
future use.

As a result, the setup algorithm only consists of two parts: the
evaluation of $E$ and $\app{W}$, and
the sparse factorization $\app{W} = L_{\app{W}}R_{\app{W}}$.
\begin{itemize}
  \item For the first part, since $E$ only has $\omega$ non-zero
    entries and computing each entry takes $O(1)$ steps
    \eqref{eq:entry}, the overall cost for this step is $O(\omega) =
    O(n)$. Once $E$ is formed, computing $\app{W}$ also takes at most
    $O(n)$ steps.
  \item The cost of the second part is more complicated. For a
    geometry that is uniformly convex, it can be shown that the number
    of non-zero entries in $L_{\app{W}}$ and $R_{\app{W}}$ is $O(n\log
    n)$. However, for a boundary with a significant flat part, the
    cost increases to $O(n^{3/2})$. The reason is that the restriction
    of $\app{W}$ to the first and last members of all $K_i$ (i.e., the
    most tangential directions) is of size $O(\sqrt{n}) \times
    O(\sqrt{n})$ and is filled significantly due to the flat
    part. Constructing LU decomposition directly for this part already
    requires $O(n^{3/2})$ steps. In order to reduce the complexity,
    the LU factorization of this submatrix is computed with the
    hierarchical matrix algebra \cite{borm-2006} or the HSS matrices
    \cite{xia-2010} as the flat part the problem is essentially a 1D
    problem. Using these hierarchical algorithms reduces the
    factorization cost to $O(n)$.
\end{itemize}
Adding these numbers together shows that the setup cost of the
preconditioner is of order $O(n)$.

Now consider the application cost of the preconditioner
\eqref{eq:eapp}. For the major steps of applying \eqref{eq:eapp}, we
  have the following estimates:
\begin{itemize}
\item The application of $\app{B^{-1}}$ is linear time due to the
  hierarchical matrix algebra approximation for $B_i^{-1}$.
\item The application of $U$ is $O(n\log n)$ since each $U_i$ is a
  partial Fourier matrix and the FFT can be used.
\item The application of $\app{T}$ is $O(n)$ since the number of
  non-zeros in $\app{T}$ is proportional to $O(n)$ after we
  threshold each $T_i$.
\item Applying $R_{\app{W}}^{-1}L_{\app{W}}^{-1}$ also takes linear
  time by using sparse backward and forward substitution algorithm,
  along with the hierarchical matrix algebra or HSS matrix for the two
  tangential submatrices.
\end{itemize}
Putting these together shows that the application cost of the
preconditioner scales like $O(n\log n)$.

%-----------------------------------
\section{Numerical Results}

\begin{figure}[h!]
  \begin{center}
    \includegraphics[height=1.5in]{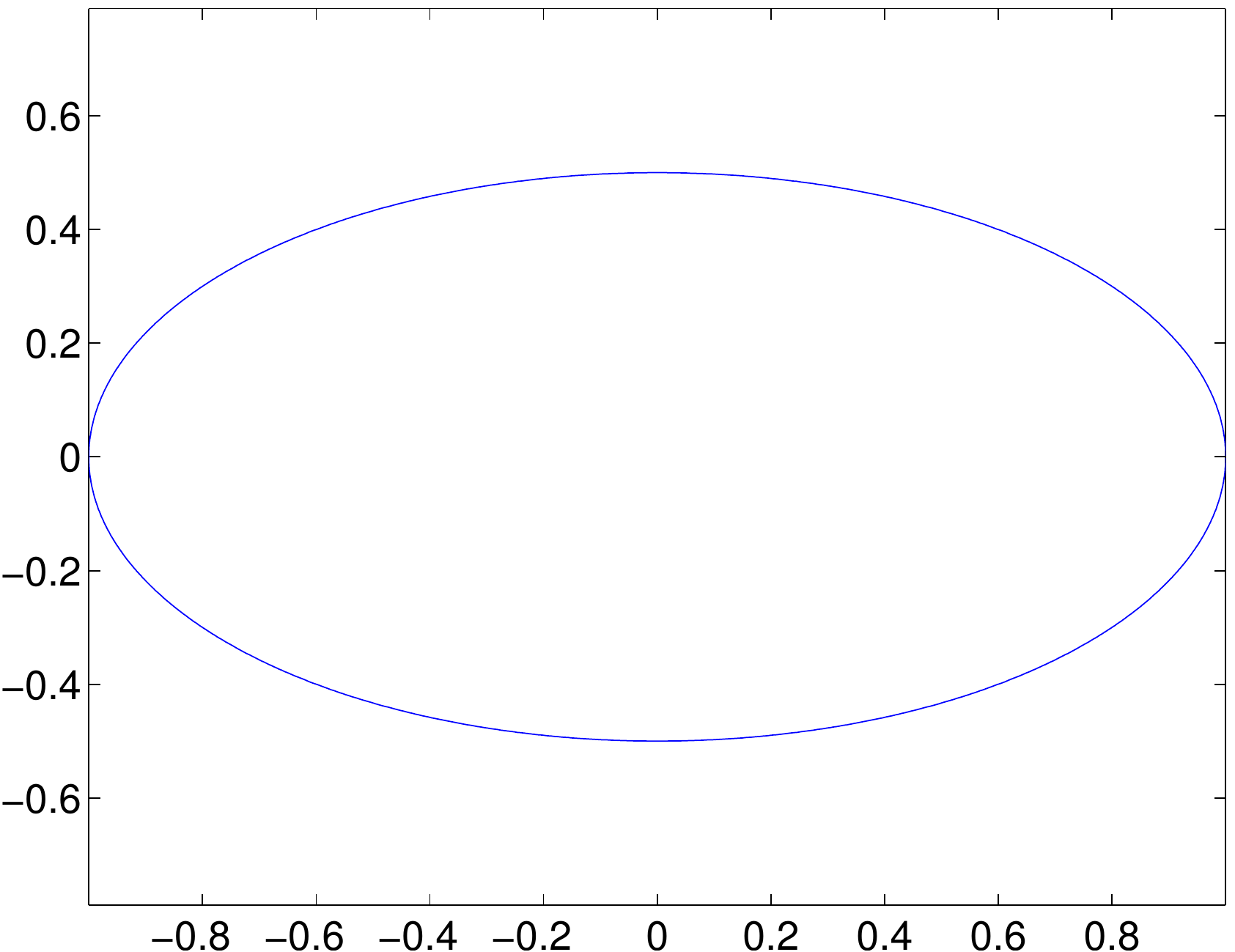}    \hspace{0.25in}
    \includegraphics[height=1.5in]{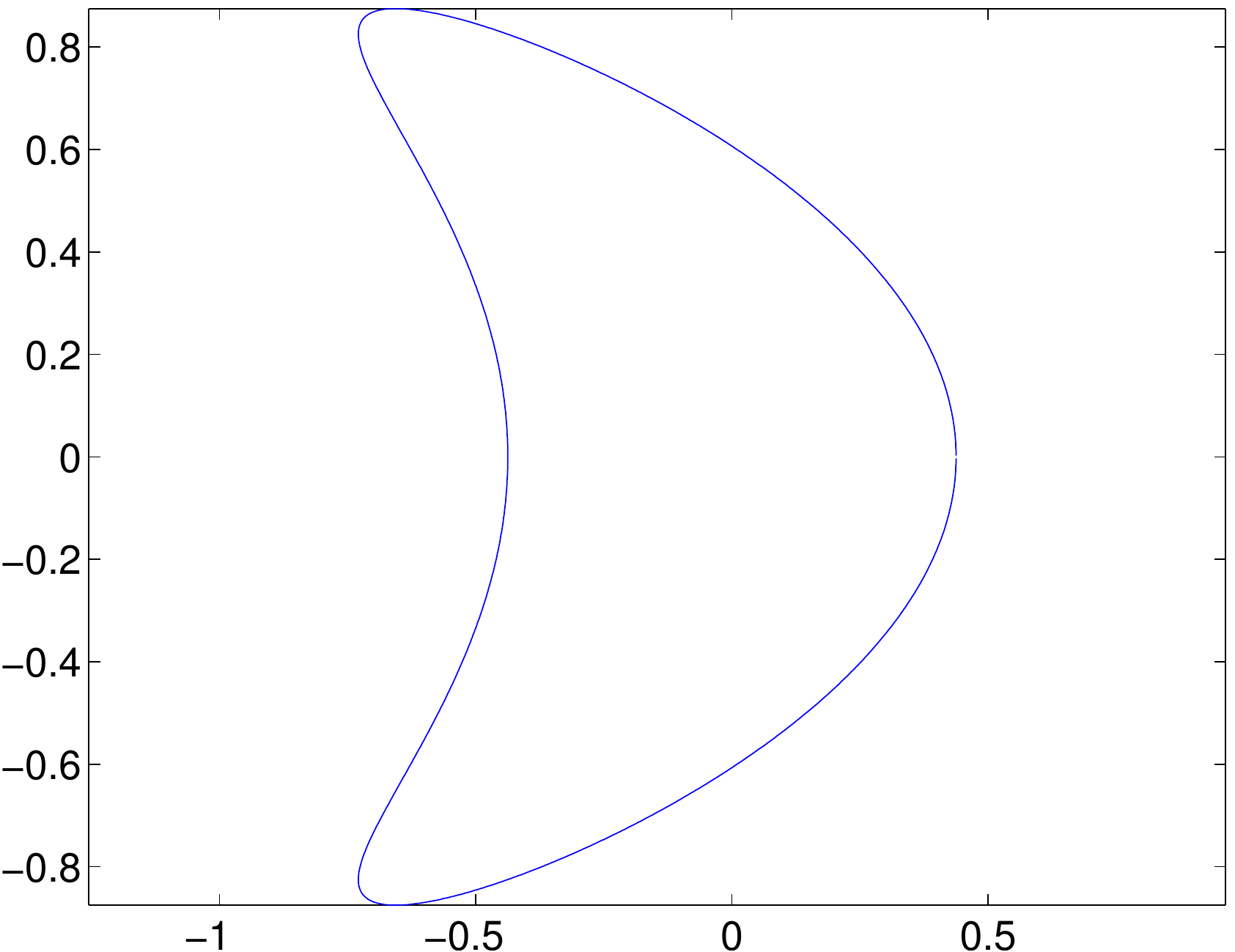}    %\hspace{0.25in}
  \end{center}
  \caption{The two scatterers used in the numerical tests. (a) an
    ellipse. (b) a bean-shaped object.}
  \label{fig:doms}
\end{figure}

The proposed preconditioner is implemented in Matlab. The numerical
results in this section are obtained on a desktop computer with a
3.60GHz CPU. Numerical tests are carried out for two domains shown in
Figure \ref{fig:doms}.

In these experiments, we discretize the domain boundary using the
Nystr\"om method with $p=8$ points per wavelength. The Chebyshev grid
size $m_c$ used to construct $e_{ij}$ is set to be $10$. The boundary
condition for both the Dirichlet and Neumann problems are obtained by
considering an incoming plane wave pointing towards the positive $x$
direction in $\R^2$.

For the iterative solution of the linear system, we use GMRES with
relative tolerance equal to $10^{-6}$ and restart number equal to
$80$. For the matrix-vector multiplication routine in the iterative
solver, the fast algorithm described in \cite{ying-2014a} is used.

We first consider the Dirichlet problem of sound-soft scattering. The
results for the two domains are reported in Tables \ref{tbl:D1} and
\ref{tbl:D2}, where
\begin{itemize}
\item $T_s$ is the setup time of the preconditioner in seconds,
\item $T_a$ is the application time of the preconditioner in seconds,
\item $T_m$ is the matrix-vector multiplication time using the
  algorithm from \cite{ying-2014a},
\item $n_p$ is the iteration number of the iterative solver when the
  preconditioner is used, and finally,
\item $n_n$ is the iteration number without any preconditioning.
\end{itemize}
The ratio $T_a/T_m$ is a good indicator of computation cost of
applying the preconditioner, while $n_p/n_n$ shows the effectiveness
of the preconditioner.

\begin{table}[ht!]
  \begin{center}
    \begin{tabular}{|cc|cc|c|cc|}
      \hline
      $\omega$ & $n$ & $T_s$ & $T_a$ & $T_m$ & $n_p$ & $n_n$\\
      \hline
      5.3e+03 & 3.3e+04 & 5.3e+00 & 2.4e-02 & 1.5e+00 & 14 & 47\\
      2.1e+04 & 1.3e+05 & 2.5e+01 & 7.4e-02 & 6.3e+00 & 16 & 71\\
      8.5e+04 & 5.2e+05 & 1.5e+02 & 3.7e-01 & 2.8e+01 & 19 & 114\\
      \hline
    \end{tabular}
  \end{center}
  \caption{Numerical results of the sound-soft scattering for the
    ellipse.}
  \label{tbl:D1}
\end{table}

\begin{table}[ht!]
  \begin{center}
    \begin{tabular}{|cc|cc|c|cc|}
      \hline
      $\omega$ & $n$ & $T_s$ & $T_a$ & $T_m$ & $n_p$ & $n_n$\\
      \hline
      5.2e+03 & 3.3e+04 & 5.3e+00 & 1.0e-02 & 1.8e+00 & 14 & 50\\
      2.1e+04 & 1.3e+05 & 2.2e+01 & 5.6e-02 & 7.7e+00 & 16 & 74\\
      8.3e+04 & 5.2e+05 & 9.2e+01 & 3.2e-01 & 3.3e+01 & 18 & 118\\
      \hline
    \end{tabular}
  \end{center}
  \caption{Numerical results of the sound-soft scattering for the
    bean-shaped object.}
  \label{tbl:D2}
\end{table}

The results show that the setup time of the preconditioner is
typically equivalent to a couple of fast matrix-vector
multiplications, while the application time of the preconditioner is
much lower. Therefore, the cost of applying the preconditioner is
almost negligible during the iterative solution. Second, the iteration
number of the preconditioned system is significantly lower than the
one of the unpreconditioned system. More importantly, the iteration
number of the former scales like $O(\log \omega)$, thus almost
frequency-independent.

Next, we consider the Neumann problem of sound-hard scattering. The
results for the two domains are reported in Tables \ref{tbl:N1} and
\ref{tbl:N2}. The results are qualitatively similar to the one for the
Dirichlet problem and demonstrate the effectiveness of the
preconditioner for the sound-hard scattering problem.

\begin{table}[ht!]
  \begin{center}
    \begin{tabular}{|cc|cc|c|cc|}
      \hline
      $\omega$ & $n$ & $T_s$ & $T_a$ & $T_m$ & $n_p$ & $n_n$\\
      \hline
      5.3e+03 & 3.3e+04 & 7.3e+00 & 1.3e-02 & 1.5e+00 & 15 & 38\\
      2.1e+04 & 1.3e+05 & 3.1e+01 & 8.2e-02 & 6.4e+00 & 19 & 56\\
      8.5e+04 & 5.2e+05 & 1.7e+02 & 4.2e-01 & 2.8e+01 & 23 & 81\\
      \hline
    \end{tabular}
  \end{center}
  \caption{Numerical results of the sound-hard scattering for the
    ellipse.}
  \label{tbl:N1}
\end{table}

\begin{table}[ht!]
  \begin{center}
    \begin{tabular}{|cc|cc|c|cc|}
      \hline
      $\omega$ & $n$ & $T_s$ & $T_a$ & $T_m$ & $n_p$ & $n_n$\\
      \hline
      5.2e+03 & 3.3e+04 & 7.3e+00 & 1.7e-02 & 1.9e+00 & 15 & 36\\
      2.1e+04 & 1.3e+05 & 2.9e+01 & 4.7e-02 & 7.7e+00 & 18 & 51\\
      8.3e+04 & 5.2e+05 & 1.2e+02 & 2.9e-01 & 3.4e+01 & 22 & 72\\
      \hline
    \end{tabular}
  \end{center}
  \caption{Numerical results of the sound-hard scattering for the
    bean-shaped object.}
  \label{tbl:N2}
\end{table}

%-----------------------------------
\section{Conclusion}

This paper presented the directional preconditioner for the combined
field integral equations (CFIEs) of high frequency acoustic obstacle
scattering in 2D. The main idea is to construct a data-sparse
approximation of the linear operator, transform it into an approximate
sparse linear system, and form an approximate inverse using efficient
sparse and hierarchical linear algebra algorithms.

We have assumed that the boundary is discretized with an equally
spaced set of discretization points. For non-equally spaced points,
the construction goes through as well, except that the FFT has to be
replaced with non-uniform FFTs. As a result, some of the
scatterer-independent precomputation can become dependent on the
discretization pattern.

A major part of future work is to extend this approach to 3D
scatterers. While the main idea should work, the lack of
equally-spaced discretization for general surfaces pose a clear
challenge for this approach.

One potential long term goal is to construct a direct solver for the
boundary integral equations of the obstacle scattering problem. It is
not clear at this point whether such a direct solver even
exists. However, this paper can be viewed a first step of exploring in
this direction.

One important ingredient of our approach is to transform a dense
oscillatory matrix to a sparse one. Once it is in a sparse form, we
can leverage the amazing power of sparse linear algebra
algorithms. While traditionally there is relatively little overlap
between the work in integral equations and the one in sparse linear
algebra, this work hints at fruitful exchange of ideas between these
two fields.

\bibliographystyle{abbrv} \bibliography{ref}

\end{document}